\documentclass[11pt]{article}
\usepackage[margin=4.46cm]{geometry}
\usepackage{graphicx}
\usepackage{epsfig}
\usepackage{amsmath}
\usepackage{amssymb}
\newcommand{\commentout}[1]{}

\def \Rset {{\mathbb R}}

\def \Nset {{\mathbb N}}

\def \Tset {{\mathbb T}}

\newcommand{\ra}{\rightarrow}

\newcommand{\nit}{\noindent}

\newcommand{\be}{\begin{equation}}
\newcommand{\ee}{\end{equation}}
\newcommand{\ba}{\begin{eqnarray}}
\newcommand{\ea}{\end{eqnarray}}
\newcommand{\bi}{\begin{itemize}}
\newcommand{\ei}{\end{itemize}}
\newcommand{\br}{\begin{eqnarray}}
\newcommand{\er}{\end{eqnarray}}

\newcommand{\dta}{\mbox{$\delta$}}

\newcommand{\eps}{\mbox{$\epsilon$}}

\newcommand{\qed}{\mbox{$\square$}\newline}

\newcommand{\norm}[1]{\lVert #1 \rVert}

\newcommand{\Rn}{\mathbb{R}^{n}}
\newcommand{\Rm}{\mathbb{R}}
\newcommand{\Rd}{\mathbb{R}^{d}}

\newtheorem{theo}{Theorem}[section]

\newtheorem{cor}{Corollary}[section]
\newtheorem{rmk}{Remark}[section]
\renewcommand{\theequation}{\arabic{section}.\arabic{equation}}

\begin{document}

\title{ A Numerical Study of   
Turbulent Flame Speeds of Curvature and Strain G-equations in Cellular Flows}
\author{Yu-Yu Liu$^{1}$, Jack Xin$^{2}$, Yifeng Yu$^{3}$ 
\thanks{$^{1,2,3}$Department of Mathematics,
UC Irvine, Irvine, CA 92697, USA. Email: yliu@math.uci.edu, jxin@math.uci.edu, yyu1@math.uci.edu.}}
\date{}

\maketitle
\thispagestyle{empty}

\begin{abstract}
We study front speeds of curvature and strain G-equations arising in turbulent combustion. 
These G-equations are Hamilton-Jacobi type level set partial differential equations (PDEs) with non-coercive Hamiltonians and degenerate nonlinear second order diffusion. 
The Hamiltonian of strain G-equation is also non-convex. Numerical computation is performed based on monotone discretization and weighted essentially nonoscillatory (WENO) approximation of transformed G-equations on a fixed periodic domain. 
The advection field in the computation is a two dimensional Hamiltonian flow consisting of a periodic array of counter-rotating vortices, or cellular flows. 
Depending on whether the evolution is predominantly in the hyperbolic or parabolic regimes, suitable explicit and semi-implicit time stepping methods are chosen. 
The turbulent flame speeds are computed as the linear growth rates of large time solutions. 
A new nonlinear parabolic PDE is proposed for the reinitialization of level set functions to prevent piling up of multiple bundles of level sets on the periodic domain. 
We found that the turbulent flame speed $s_T$ of the curvature G-equation is enhanced as the intensity $A$ of cellular flows increases, at a rate between those of the inviscid and viscous G-equations. 
The $s_T$ of the strain G-equation increases in small $A$, decreases in larger $A$, then drops down to zero at a large enough but finite value $A_{*}$. 
The flame front ceases to propagate at this critical intensity $A_*$, and is quenched by the cellular flow.       
\end{abstract}

\hspace{.1 in} {\bf Key Words}: 
Curvature/Strain G-equations, Cellular Flows, 

\hspace{.1 in} 
Front Speed Computation, Enhancement and Quenching. 
\medskip

\hspace{.1 in} {\bf AMS Subject Classification:} 
70H20, 76F25, 76M50, 76M20.

\newpage

\section{Introduction}
\setcounter{equation}{0}
Front propagation in turbulent combustion is a nonlinear and multiscale dynamical process \cite{W85,Yak,S89,Siv,MS_94,MS_98,P00,Xin_00}. 
The first principle based approach requires a system of reaction-diffusion-advection equations coupled with the Navier-Stokes equations. 
Simplified models, such as the advective Hamilton-Jacobi equations (HJ) and passive scalar reaction-diffusion-advection
equations (RDA), are often more efficient in improving our understanding of such complex phenomena. 
Progress is well documented in books \cite{W85,P00,Xin_09} and research papers \cite{Ab_02,ABP_00,CW,Const,EMS,MS_98,NX_09a,OF02,R95,Siv,S89,Xin_00,Yak} among others.
\medskip

A sound phenomenological approach in turbulent combustion is the level set formulation \cite{OF02} of flame front motion laws with the front width ignored \cite{P00}.
The simplest motion law is that the normal velocity of the front ($V_n$) is equal to a constant $s_L$ (the laminar speed) plus the projection of fluid velocity $V(x,t)$ along the normal $\vec{n}$. 
The laminar speed is the flame speed due to chemistry (reaction-diffusion) when fluid is at rest.
As the fluid is in motion, the flame front will be wrinkled by the fluid velocity.
Under suitable conditions, the front location eventually moves to leading order at a well-defined steady speed $s_T$ in each specified direction, which is the so-called "turbulent burning velocity" \cite{P00}. 
The study of existence and properties of turbulent flame speed $s_T$ is a fundamental problem in turbulent combustion theory and experiments \cite{W85,R95,P00}. 
Let the flame front be the zero level set of a function $G(x,t)$, then the normal direction is $DG/|DG|$ and the normal velocity is $-G_t/|DG|$. ($D$: spatial gradient.) 
The motion law becomes the so-called $G$-equation in turbulent combustion \cite{W85,P00}:
\be
G_t + V(x,t)\cdot DG + s_L |DG|=0. \label{Gi}
\ee
Chemical kinetics and diffusion rates are all included in the laminar speed $s_L$ which is provided by a modeler. 
Formally under the G-equation model, for a specified unit direction $P$,
\be
s_T(P)=-\lim_{t\to +\infty}{G(x,t)\over t}. \label{speed}
\ee
Here $G(x,t)$ is the solution of equation (\ref{Gi}) with initial data $G(x,0)=P\cdot x$. 
The existence of $s_T$ has been rigorously established in \cite{XY_10} and \cite{CNS} independently for incompressible periodic flows, and \cite{NN_11} for two dimensional incompressible random flows.
\medskip

As fluid turbulence is known to cause stretching and corrugation of flames, additional modeling terms may be incorporated into 
the basic G-equation (\ref{Gi}). 
In this paper, we shall study turbulent burning velocity $s_T$ of such extended G-equation models involving strain and curvature effects. 
The curvature G-equation is:
\be\tag{Gc}\label{Gc}
G_t+V(x,t)\cdot DG+s_L|DG|
=ds_L|DG|\mathrm{div}\left({DG\over|DG|}\right),
\ee
which comes from adding mean curvature term to the basic motion law. 
The curvature dependent motion is well-known, see \cite{OS88,OF02} and references therein.
If the curvature term is further linearized \cite{Denet99}, 
we arrive at the viscous G-equation:
\be\tag{Gv}\label{Gv}
G_t+V(x,t)\cdot DG+s_L|DG|
=ds_L\Delta G,
\ee
which is also a model for understanding numerical diffusion \cite{OF02}.
The strain G-equation is:
\be\tag{Gs}\label{Gs}
G_t+V(x,t)\cdot DG+\left(s_L+
d{DG\cdot DV\cdot DG\over |DG|^2}\right)|DG|
=ds_L|DG|\mathrm{div}\left({DG\over|DG|}\right).
\ee
The strain term $n\cdot DV\cdot n$ will be derived and analyzed later. 
The formula (\ref{speed}) formally extends to (\ref{Gc}) and (\ref{Gs}). 
A complete mathematical theory of their existence is lacking at the moment. 
Helpful empirical observations from experiments \cite{B92,R95} are: 
(i) When the flame front is wrinkled by the advection, the interface area increases and $s_T$ increases (called "enhancement").
(ii) However, turbulent flame speed cannot increase without limit, and the growth rate may be sublinear in the large intensity limit of the advection (called "bending"). 
(iii) When the advection is strong up to certain level, the reactant totally scatters.  
The reaction then fails and the flame front extinguishes (called "quenching"). 
\medskip

We aim to understand and quantify these nonlinear phenomena in the context of curvature and strain G-equations and cellular flows where $s_T$ is related to the corrector (cell) problem of 
homogenization theory for which several mathematical results are available. The cellular flow is a two dimensional incompressible flow:
\be\label{Cell} 
V=\nabla^{\perp}\mathcal{H}=(-\mathcal{H}_y,\mathcal{H}_x)
\ ,\ \mathcal{H}={A\over 2\pi}\sin(2\pi x)\sin(2\pi y),
\ee
where $A$ is the amplitude of the flow. 
By parameterizing $s_T$ as a function of $A$, we are interested in the behavior of $s_T$ as $A$ increases in G-equations (\ref{Gi}),(\ref{Gc}),(\ref{Gv}),(\ref{Gs}).
The streamlines of the cellular flow consist of a periodic array of hyperbolic (saddle points, separatrices) and elliptic (vortical) regions.  
For inviscid G-equation (\ref{Gi}), it is known \cite{O02,Ab_02,NXY} that $s_T=O({A/\log(A)})$, where the logarithmic factor is due to slow-down of transport near saddle points. 
For viscous G-equation, we recently proved \cite{LXY11} that $s_T=O(1)$ as $A \gg 1$ at any fixed positive viscosity ($d > 0$). 
The dramatic slowdown (strong bending) is due to the smoothing of the level set function $G$ by viscosity, and the uniform bound of $\| DG\|_{L^{1}_{\rm loc}} $. 
Less is known about the growth rate of $s_T$ for curvature and strain G-equations. 
The curvature term only provides partial smoothing, hence the slowdown (bending) is weaker in general than the regular smoothing by viscosity. 
For shear flows, we showed \cite{LXY11} that the linear growth rate $\lim_{A\ra \infty} s_T/A$ is same as that of the inviscid G-equation. 
The effect of strain term is more difficult to analyze, as it is 
highly nonlinear in $G$ and can take both signs. 
It also changes the type of the Hamiltonian of G-equation from convex in (\ref{Gi}) to non-convex in (\ref{Gs}). 
For shear flows, the strain term always slows down $s_T$ \cite{XY_11}. 
\medskip

We shall first approximate the G-equations by a monotone discrete system, then apply high resolution numerical methods such as WENO (weighted essentially non-oscillatory finite 
difference methods \cite{JP00, OF02}) with a combination of explicit and semi-implicit time stepping strategies, 
depending on the size and property of dissipation in the equations. 
The computation is done on transformed G-equations over a periodic domain to avoid the need of excessively large computational domains to contain potentially fast moving fronts. 
We also devise a new reinitialization equation on the periodic domain to prevent the level sets from piling up during time evolution. 
A nonlinear diffusion term is added to the standard reinitialization equation (Chapter 7, \cite{OF02}) to perform reinitialization on multiple bundles of level sets often encountered during long time computation. 
An iterative method of computing $s_T$ of the viscous G-equation (\ref{Gv}) works well based on the corrector equation of homogenization, if the viscosity $d$ is above a certain level. 
\medskip

Our main findings are: 
(1) The curvature G-equation (\ref{Gc}) always enhances $s_T$ as $A$ increases; the amount of enhancement is smaller than that of the inviscid G-equation (\ref{Gi}), larger than that of the viscous G-equation (\ref{Gv}). 
For small enough $d$, the $s_T$ of (\ref{Gc}) behaves similarly to that of the inviscid G-equation (\ref{Gi}), or weak speed bending. 
For large enough $d$, the $s_T$ of (\ref{Gc}) behaves similarly to that of the viscous G-equation (\ref{Gv}), or strong speed bending. 
(2) The $s_T$ is a monotone decreasing function of $d$ for both curvature and strain G-equations, (\ref{Gc}) and (\ref{Gs}). 
(3) For the strain G-equation (\ref{Gs}) with fixed $d >0$, $s_T$ first increases with $A$, then decreases in $A$, and drops down to zero at finite $A$ (front quenching occurs). 
\medskip

The paper is organized as follows. 
In section 2, we give a brief derivation of G-equation models and an an overview of analytical results of the turbulent flame speeds. 
In section 3, we introduce numerical scheme for each G-equation. 
We also discuss how to perform reinitialization in periodic domain. 
In section 4, we present and interpret the numerical results. 
Concluding remarks are in section 5. 
In the two appendices, we show a formula of surface stretch rate in advection and a convergent iteration scheme of $s_T$ based on the corrector problem of homogenization.
\medskip

The work was partially supported by NSF grants DMS-0911277 (JX) and DMS-0901460 (YY). 
YL thanks the Department of Mathematics of UC Irvine for a  graduate fellowship. 
The authors would like to thank Professor John Lowengrub for helpful conversations on interface computations.  

\section{Derivation and Analysis of G-equations}
\setcounter{equation}{0}
\subsection{G-equations}
In the thin reaction zone regime and the corrugated flamelet regime of premixed turbulent combustion (pp. 91-107, Chapter 2, \cite{P00}), the flame front is modeled by a level set function:
$\{(x,t): G(x,t)=0 \}$, which is the interface between the burned area $\{G<0\}$ and the unburned area $\{G>0\}$. 
See \cite{OF02} for an introduction on level set methods in a broad context.   
The unit normal direction is $n=DG/|DG|$ and the normal velocity is $-G_t/|DG|$. ($D$: spatial gradient.) 
The simplest motion law is that the normal velocity of the interface is the sum of a constant $s_L$ (called {\it laminar flame speed}) and the projection of fluid velocity $V(x,t)$ 
along the normal direction. 
The $s_L$ is well-defined if the reaction zone is much larger 
than the smallest turbulent length scale (the Kolmogorov scale),
as in the corrugated flamelet regime \cite{P00}. 
In terms of $G$, the law is the so-called $G$-equation (\ref{Gi}). 
A linear version dated back to \cite{M64}. 
The trajectory of a particle $x(t)$ on the interface satisfies: 
\be \label{Ge2}
{dx\over dt}=V(x,t)+s_Ln.
\ee

The G-equation or level set framework is a popular and robust phenomenological approach. 
The motion law is in the hands of a modeler based on theory and experiments. 
Various nonlinear effects may be built into the basic model (\ref{Gi}). 
For example, turbulence is known to cause stretching of flame fronts. 
It was shown in \cite{MM82,M83} that the flame stretch rate may be added as a first order correction term on the laminar flame speed:
\be\label{Ge3}
\hat s_L=s_L-d{1\over\sigma}{d\sigma\over dt},
\ee
where $\sigma$ is the surface element area of the level set and 
$d$ is called the {\it Markstein diffusive number}. 
If the flame stretch rate is positive, the reactant on the flame front scatters and the burning reaction slows down. 
By a kinematic calculation (see appendix A for details), the flame stretch rate is:
\be\label{Ge4}
{1\over\sigma}{d\sigma\over dt}=\mathcal{S}+s_L\kappa,
\ee
$$
\mathcal{S}=-n\cdot DV\cdot n\ ,\ \kappa=\mathrm{div}(n),
$$
where $\mathcal{S}$ is called the {\bf strain rate} and $\kappa$ is the mean curvature of level set. 
Replacing $s_L$ by $\hat{s}_L$ in (\ref{Ge2}), we have the {\bf strain G-equation} (\ref{Gs}). 
In the thin reaction zone regime (section 2.6, pp 104--107 \cite{P00}), Kolmogorov scale eddies enter the reaction zone, and cause unsteady perturbations of laminar speed $s_L$. 
The $(s_L - d\mathcal{S})$ term and the eddy effects are lumped together as a fluctuating quantity (denoted by $s_{L,s}$ in \cite{P00}) which however is on the order of $s_L$ based on direct numerical simulation data.  
If we approximate $s_{L,s}$ by $s_L$ and keep the curvature term, the {\bf curvature G-equation} (\ref{Gc}) follows. 

\begin{rmk}
In previous works \cite{ZR94,AS97}, $\hat{s}_L$ is modified to remain positive:
$$
\hat{s}_L=\max\left\{s_L-d{1\over\sigma}{d\sigma\over dt},0\right\}\ ,\ 
s_L\exp\left(-{d\over s_L}{1\over\sigma}{d\sigma\over dt}\right).
$$
However, these modifications restrict the curvature or strain effect in the strong advection scheme. 
The bending or quenching effect may either weaken or disappear.
\end{rmk}


\subsection{Turbulent Burning Velocity}
We discuss how to evaluate turbulent flame speeds in G-equation models. 
For simplicity we consider the inviscid G-equation (\ref{Gi}) only, and the formulation extends to other G-equations. 

Given a unit vector $P\in\Rn$ and suppose the flame front propagates in direction $P$.
Let the initial flame front be $\{P\cdot x=0\}$ and consider G-equation with planar initial condition:
\be\label{Tu0}
\left\{\begin{array}{ll}
G_t +V(x,t)\cdot DG+ s_L|DG| = 0 & 
\mbox{in}\ \mathbb{R}^n\times (0,\infty)
\vspace{0.05in}\\
G(x,0) = P\cdot x &  
\mbox{on}\ \mathbb{R}^n\times \{t=0\}
\end{array}. \right.
\ee
Assume $V(x,t)$ is spatially periodic. If we write $G(x,t)= P\cdot x + u(x,t)$, then $u(x,t)$ is also spatially periodic and solves the following periodic initial value problem:
\be\label{Tu3}
\left\{\begin{array}{ll}
u_t +V(x,t)\cdot (P+Du)+ s_L|P+Du| = 0 & 
\mbox{in}\ \mathbb{T}^n\times (0,\infty)
\vspace{0.05in}\\
u(x,0) = 0 &  
\mbox{on}\ \mathbb{T}^n\times \{t=0\}
\end{array}.\right.
\ee
Hence in numerical computation of (\ref{Tu0}) we can reduce the spatial domain from $\Rn$ to $[0,1]^n$ by imposing the affine periodic condition:
\be\label{Tu1}
\left\{\begin{array}{ll}
G_t +V(x,t)\cdot DG+ s_L|DG| = 0 & 
\mbox{in}\ [0,1]^n\times (0,\infty)
\vspace{0.05in}\\
G(x,0) = P\cdot x &  
\mbox{on}\ [0,1]^n\times \{t=0\}
\end{array}. \right.
\ee
\be\label{Tu2}
G(x+z,t) = G(x,t)+P\cdot z,\ x\in [0,1]^n,\ z\in\mathbb{Z}^n.
\ee

Now we focus on $P=e_1=(1,0,\cdots, 0)$, then $G(x,t)=x_1+u(x,t)$ is periodic in $x_2,\cdots, x_n$. 
Consider the stripe domain $\mathbb{R}\times [0,1]^{n-1}$, and the burned area at time $t$ is $\{x\in \mathbb{R}\times [0,1]^{n-1}:G(x,t)<0\}$. 
Denote $\mathcal{A}(t)$ the volume that burned area has invaded during time interval $(0,t)$, then turbulent flame speed is the linear growth rate of $\mathcal{A}(t)$:
\be\label{Tu4}
s_T=\lim_{t\ra+\infty}{\mathcal{A}(t)\over t}
=\lim_{t\ra+\infty}{1\over t}\int_{\mathbb{R}\times [0,1]^{n-1}}
\left(\chi_{\{G(x,t)<0\}}-\chi_{\{G(x,0)<0\}}\right)dx.
\ee 
($\chi$: indicator function.) 
Note that $G(x,0)=x_1$ and $G(x+e_1,t)=G(x,t)+1$, then $\mathcal{A}(t)$ and hence $s_T$ can be evaluated by $G$ or $u$ in $[0,1]^n$:
\be\label{Tu5}
s_T
=\lim_{t\ra+\infty}
{-1\over t}\int_{[0,1]^n}\left[G(x,t)\right]dx
=\lim_{t\ra+\infty}
{-1\over t}\int_{[0,1]^n}\left[x_1+u(x,t)\right]dx.
\ee
($[\cdot]$: floor function.) 
In \cite{ZR94} the initial condition is chosen as $G(x,0)=\phi(x_1)$ with $\phi:\Rset\ra\Rset$ a smeared-out signed function, and the computational domain is 
$[a,b]\times [0,1]$. 
If the zero level set travels a long distance, the length of the domain ($b-a$) needs to be large enough to contain the level set. 
To study a fast moving flame front and its long time behavior, 
the computational domain will be very large. 
Instead we choose $G(x,0)=x_1$ and reduce the computational domain to $[0,1]\times [0,1]$. The $s_T$ is the same from either initial data.


Another way to find turbulent flame speed is via the framework of periodic homogenization \cite{LPV88, E92}. 
Assume $V=V(x)$ be time-independent periodic flow and consider the so-called corrector problem: given any vector $P\in\Rn$, 
find a number $\bar{H}$ (the effective Hamiltonian) such that the equation
\be\label{Tu6}
V(x)\cdot (P+D\bar{u})+ s_L|P+D\bar{u}| = \bar{H},\ x\in\mathbb{T}^n
\ee
has a periodic solution $\bar{u}(x)$. 
If (\ref{Tu6}) is solvable, then  G-equation has the following stationary solution: 
\be\label{Tu7}
G(x,t)=-\bar{H}t+ P\cdot x + \bar{u}(x),
\ee 
and $\bar{H}$ is exactly the turbulent flame speed. 
The corrector problem is well-posed for viscous G-equation \cite{LXY10}, and can be used to compute $s_T$ iteratively when viscosity is not too small (see section 3.4 and Appendix B).  
However, (\ref{Tu6}) for inviscid G-equation may not have exact solutions due to lack of coercivity of G-equations, only approximate solutions exist \cite{XY_10}. 
It is also an open question in general whether it has solutions if the curvature or strain term is present. 
The more general and robust characterization of $s_T$ is simply the linear growth rate of $G$ or $u$ at fixed $x$:
\be\label{Tu8}
s_T=\lim_{t\ra+\infty}{-G(x,t)\over t}=\lim_{t\ra+\infty}{-u(x,t)\over t},
\ee
which we shall adopt for curvature and strain G-equations in this paper. 
Indeed (\ref{Tu5}) and (\ref{Tu8}) are consistent when $P=e_1$, but (\ref{Tu8}) can be used for any direction $P$.
See \cite{Q03} for earlier work on computing effective Hamiltonian of coercive Hamilton-Jacobi equations along this line. 

\section{Numerical Methods}
\setcounter{equation}{0}
We discuss the numerical schemes for G-equations. 
We employ the Hamilton-Jacobi weighted essentially nonoscillatory (HJ WENO) scheme and the total variation diminishing Runge-Kutta (TVD RK) scheme in higher order spatial and time discretization respectively. 
See \cite{JP00, SO88, OF02} for details of the schemes.

\subsection{Inviscid G-equation}
Inviscid G-equation (\ref{Gi}) is a Hamilton-Jacobi equation with Hamiltonian
\be\label{Inv1}
H(p)=V(x,t)\cdot p+s_L|p|.
\ee
The forward Euler time discretization of (\ref{Gi}) is
\be\label{Inv2}
{G^{n+1}-G^n\over\Delta t}+
\hat{H}^n(G_x^-,G_x^+,G_y^-,G_y^+)=0,
\ee
where $\hat{H}$ is the numerical Hamiltonian of (\ref{Inv1}) and $G_x^-$ ($G_x^+$) denotes the left (right) discretization of $G_x$. 
For the stability of the numerical scheme, $\hat{H}=\hat{H}(p_x^-,p_x^+,p_y^-,p_y^+)$ is chosen to be consistent and monotone \cite{CL84}. 
Here consistency means that $\hat{H}(p_x,p_x,p_y,p_y)=H(p_x,p_y)$, and monotonicity means that $\hat{H}$ is nondecreasing in $p_x^-,p_y^-$ and nonincreasing in $p_x^+,p_y^+$.

Write $V=(V_1,V_2)$ in (\ref{Inv1}):
$$
H(p_x,p_y)=\left(V_1+s_L{p_x\over|p|}\right)p_x+
\left(V_2+s_L{p_y\over|p|}\right)p_y.
$$
When the velocity field dominates the normal velocity, upwinding direction is determined by the velocity field. 
For example, if $V_1>s_L$, then $V_1+s_Lp_x/|p|$ is always positive and $p_x$ is approximated by $p^-_x$. 
However, if the velocity field and the normal velocity are comparable, it is hard to determine the upwinding direction.
In this case we treat both terms separately: for the velocity field term, we apply upwinding scheme; for the normal velocity term, we apply Godunov scheme. 
Since both schemes are monotone, their sum is again monotone.
In summary, we have the following monotone numerical Hamiltonian of (\ref{Inv1}):
\be\label{Inv3}
\hat{H}(p_x^-,p_x^+,p_y^-,p_y^+)
= V_1p_x^{Vel}+V_2p_y^{Vel}
+s_L\sqrt{(p_x^{Nor})^2+(p_y^{Nor})^2},
\ee
where
$$
p_x^{Vel}=\left\{ \begin{array}{ll}
p_x^- &,\mbox{if}\ V_1>0\vspace{0.05in}\\
p_x^+ &,\mbox{if}\ V_1<0
\end{array}\right.\ , \
p_y^{Vel}=\left\{ \begin{array}{ll}
p_y^- &,\mbox{if}\ V_2>0\vspace{0.05in}\\
p_y^+ &,\mbox{if}\ V_2<0
\end{array}\right.
$$
and
$$
\begin{array}{l}
(p_x^{Nor})^2=\left\{ \begin{array}{ll}
(p_x^-)^2 &,\mbox{if}\ V_1>s_L\vspace{0.05in}\\
\max\left(\max(p_x^-,0)^2,\min(p_x^+,0)^2\right) &,\mbox{if}\ |V_1|\leq s_L\vspace{0.05in}\\
(p_x^+)^2 &,\mbox{if}\ V_1<-s_L
\end{array}\right.,\vspace{0.05in}\\
(p_y^{Nor})^2=\left\{ \begin{array}{ll}
(p_y^-)^2 &,\mbox{if}\ V_2>s_L\vspace{0.05in}\\
\max\left(\max(p_y^-,0)^2,\min(p_y^+,0)^2\right) &,\mbox{if}\ |V_2|\leq s_L\vspace{0.05in}\\
(p_y^+)^2 &,\mbox{if}\ V_2<-s_L
\end{array}\right..
\end{array}
$$

For the accuracy of the numerical scheme, we apply WENO5 scheme to approximate the spatial derivatives and RK3 scheme in forward Euler time discretization. 
The time step restriction (CFL condition) is
\be\label{Inv4}
\Delta t\left(
{\|V_1\|+s_L\over \Delta x}+{\|V_2\|+s_L\over \Delta y}
\right)<1.
\ee
($\|\cdot\|$: maximum norm.) Overall the scheme gives nearly fifth order spacial accuracy in smooth regions of solutions, and third order accuracy in time.

\begin{rmk}
Compared with the standard schemes (LF, LLF, RF, etc. See chapter 5 of \cite{OF02}), our choice of numerical Hamilton is easy to implement, and no extra artificial diffusion is added to satisfy the monotonicity.
\end{rmk}

\subsection{Curvature G-equation}
In forward Euler scheme of curvature G-equation (\ref{Gc}), the curvature term in two dimensional space is
$$
|DG|\mathrm{div}\left({DG\over|DG|}\right)
={G_y^2G_{xx}-2G_xG_yG_{xy}+G_x^2G_{yy}\over G_x^2+G_y^2}
$$
and is discretized by central differencing \cite{OS88}.
Since central differencing gives only second order accuracy, we apply WENO3 scheme to evaluate the numerical Hamiltonian (\ref{Inv3}) and RK2 scheme in time step discretization. 
The time step restriction is
\be\label{Cur1}
\Delta t\left(
{\|V_1\|+s_L\over \Delta x}+{\|V_2\|+s_L\over \Delta y}
+{2s_Ld\over(\Delta x)^2}+{2s_Ld\over(\Delta y)^2}
\right)<1.
\ee

When $d$ is large ($\gg \Delta x$), the time step size for forward Euler scheme is very small $\Delta t=O((\Delta x)^2)$. To alleviate the stringent time step restriction, we decompose the curvature term as follows:
\be\label{Cu2}
|DG|\mathrm{div}\left({DG\over|DG|}\right)
=\Delta G-\Delta_{\infty} G
\ee
$$
=(G_{xx}+G_{yy})-{G_x^2G_{xx}+2G_xG_yG_{xy}+G_y^2G_{yy}\over G_x^2+G_y^2},
$$
where $\Delta_{\infty}$ is the infinity Laplacian operator. If we apply backward Euler scheme on $\Delta G$ and forward Euler scheme on $\Delta_{\infty} G$, then we have the following semi-implicit time discretization scheme for (\ref{Gc}):
\be\label{Cu3}
{G^{n+1}-G^n\over \Delta t}+V(x,t^n)\cdot DG^n+s_L|DG^n|
=ds_L(\Delta G^{n+1}-\Delta_{\infty} G^n),
\ee
whose time step restriction is same as inviscid G-equation (\ref{Inv4}). 
Note that for implicit scheme each time step is more expensive. Hence if $d$ is small ($\sim \Delta x$), the forward Euler scheme is still the better choice.

Another cause of small time step is when $\|V\|$ is large.
However we cannot move the advection term into implicit scheme as in standard advection-diffusion equations. 
The curvature G-equation is essentially of hyperbolic type rather than of parabolic type. 
Even involving second order derivatives, the curvature term is dissipative only along the tangential plane of the level set and so cannot stabilize the advection term.
 
\begin{rmk} The curvature term and the infinity Laplacian operator in higher dimensional space are
$$
|DG|\mathrm{div}\left({DG\over|DG|}\right)
=\left(\dta_{ij}-{G_{x_i}G_{x_j}\over|DG|^2}\right)G_{x_{i}x_{j}}\ ,\ 
\Delta_{\infty}G={G_{x_i}G_{x_j}\over|DG|^2}G_{x_{i}x_{j}}.
$$
\end{rmk}

\subsection{Strain G-equation}
For strain G-equation (\ref{Gs}), the Hamiltonian becomes
\be\label{Str1}
H(p)=V(x,t)\cdot p+(s_L-d\mathcal{S})|p|\ ,\ 
\mathcal{S}=-{p\cdot DV\cdot p\over |p|^2}.
\ee
If we apply upwinding scheme on $V\cdot p$, then it suffices to find a monotone scheme for $(s_L-d\mathcal{S})|p|$.
First we approximate $p$ to obtain $\mathcal{S}$, and next we evaluate $|p|$ by Godunov scheme according to the sign of $(s_L-d\mathcal{S})$. 
Then we obtain the following monotone numerical Hamiltonian of (\ref{Str1}):
\be\label{Str2}
\hat{H}(p_x^-,p_x^+,p_y^-,p_y^+)
= V_1p_x^{Vel}+V_2p_y^{Vel}
+(s_L-d\hat{\mathcal{S}})\sqrt{(p_x^{Nor})^2+(p_y^{Nor})^2},
\ee
where $p_x^{Vel}$, $p_y^{Vel}$ are same as in (\ref{Inv3}), $\hat{\mathcal{S}}$ is the numerical approximation of $\mathcal{S}$ with $p$ evaluated by central differencing, and
$$
(p_x^{Nor})^2=
\left\{\begin{array}{l}
\max(\max(p_x^-,0)^2,\min(p_x^+,0)^2),
\ \ \mbox{if}\ \ (s_L-d\hat{\mathcal{S}})>0
\vspace{0.05in}\\
\max(\min(p_x^-,0)^2,\max(p_x^+,0)^2),
\ \ \mbox{if}\ \ (s_L-d\hat{\mathcal{S}})<0
\end{array}\right.,
$$
$$
(p_y^{Nor})^2=
\left\{\begin{array}{l}
\max(\max(p_y^-,0)^2,\min(p_y^+,0)^2),
\ \ \mbox{if}\ \ (s_L-d\hat{\mathcal{S}})>0
\vspace{0.05in}\\
\max(\min(p_y^-,0)^2,\max(p_y^+,0)^2),
\ \ \mbox{if}\ \ (s_L-d\hat{\mathcal{S}})<0
\end{array}\right..
$$

\begin{rmk}
For cellular flow (\ref{Cell}), the strain rate can be simplified as
$$
\mathcal{S}=-2\pi A\cos(2\pi x)\cos(2\pi y)
{(G_y^2-G_x^2)\over |DG|^2}.
$$
Then $(s_L-d\mathcal{S})$ is always positive if $2\pi Ad<s_L$.
\end{rmk}


\subsection{Viscous G-equation}
When $d$ is small, viscous G-equation (\ref{Gv}) is advection dominated and should be treated like a hyperbolic equation. 
Similar to curvature G-equation, for spatial discretization, we apply WENO3 scheme on numerical Hamiltonian (\ref{Inv3}) and central differencing on the diffusion term. 
For time step discretization, we apply RK2 forward Euler scheme.

When $d$ is large enough, we consider the following semi-implicit scheme:
\be\label{Vis1}
{G^{n+1}-G^{n}\over\Delta t}+V(x,t^{n+1})\cdot DG^{n+1}+s_L|DG^n|
=ds_L\Delta G^{n+1}.
\ee
Here the advection and diffusion terms are discretized by central differencing, and the normal direction term is discretized by Godunov and WENO3 scheme. 
Since there is no time step restriction from both advection and diffusion terms, the time step constraint for (\ref{Vis1}) is
$$
\Delta t\left(
{s_L\over \Delta x}+{s_L\over \Delta y}
\right)<1.
$$

When $V=V(x)$ is periodic, mean zero and incompressible, turbulent flame speed may also be obtained from the corrector problem:
\be\label{Vis2}
-ds_L\Delta \bar{u}+V(x)\cdot (P+D\bar{u})+ s_L|P+D\bar{u}| = \bar{H}
,\ x\in\Tset^n, 
\ee
which has a unique (up to a constant) classical solution and $\bar{H}=s_L\int_{\Tset^n}|P+D\bar{u}|dx$. 
When $d$ is large enough, the following iteration scheme converges:
$$
-ds_L\Delta u^{(k+1)}+V(x)\!\cdot\! Du^{(k+1)}
= H^{(k)}-s_L|P+Du^{(k)}|-V(x)\!\cdot\! P,\ x\in\Tset^n,
$$
\be \label{Vis3}
H^{(k)}=s_L\int_{\Tset^n}|P+Du^{(k)}|dx.
\ee
A convergence proof is in Appendix B. 
To solve (\ref{Vis3}) numerically as an elliptic equation, all operators are discretized by central differencing.


\subsection{Reinitialization}

\begin{figure}
\center (a) \\
\includegraphics[width=0.8\textwidth]{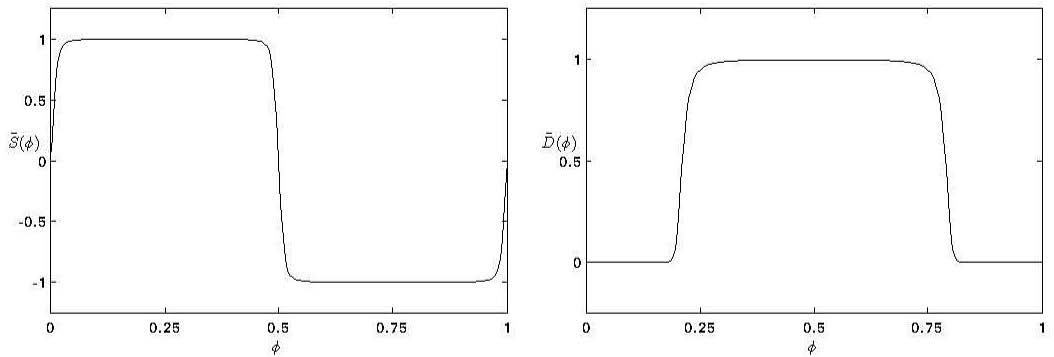}
\center (b) \hspace{2in} (c)\\
\includegraphics[width=0.48\textwidth]{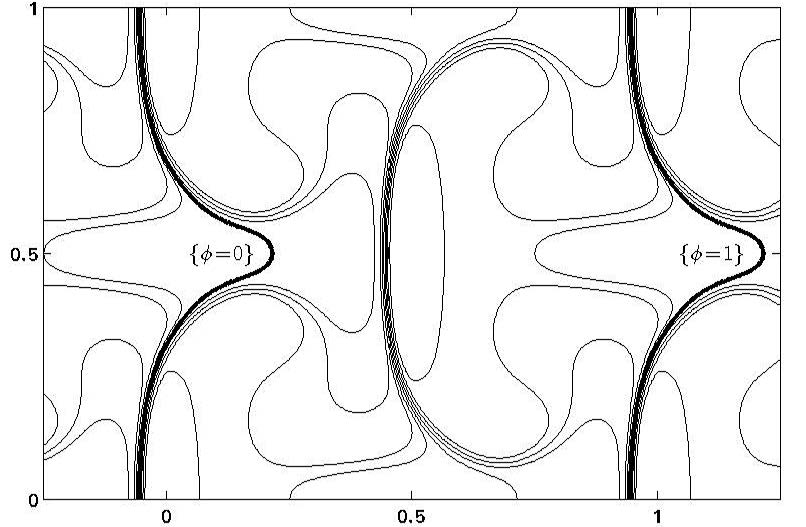}
\includegraphics[width=0.48\textwidth]{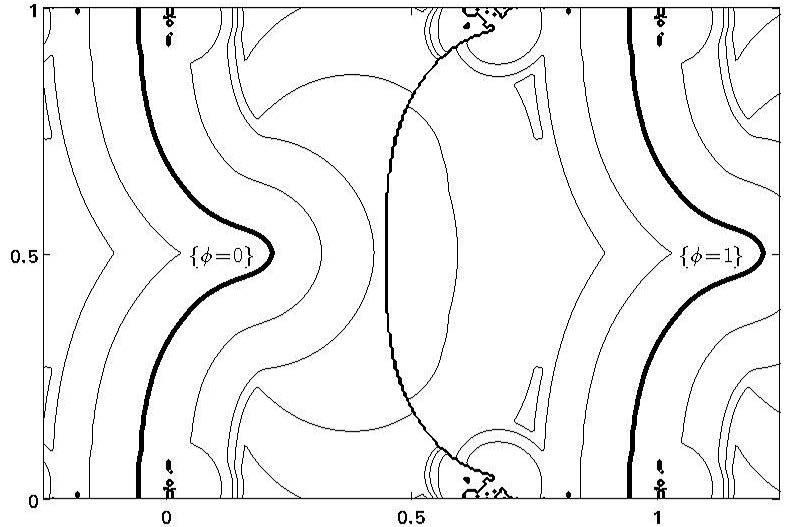}
\center (d) \hspace{2in} (e)\\
\includegraphics[width=0.48\textwidth]{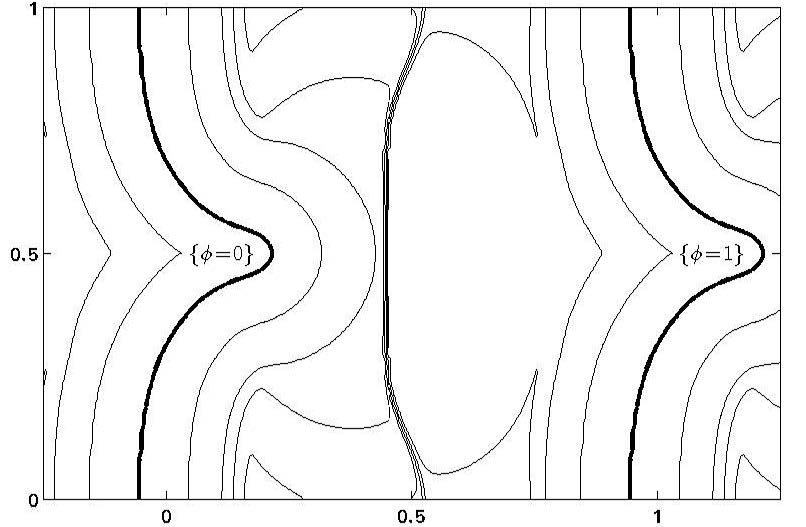}
\includegraphics[width=0.48\textwidth]{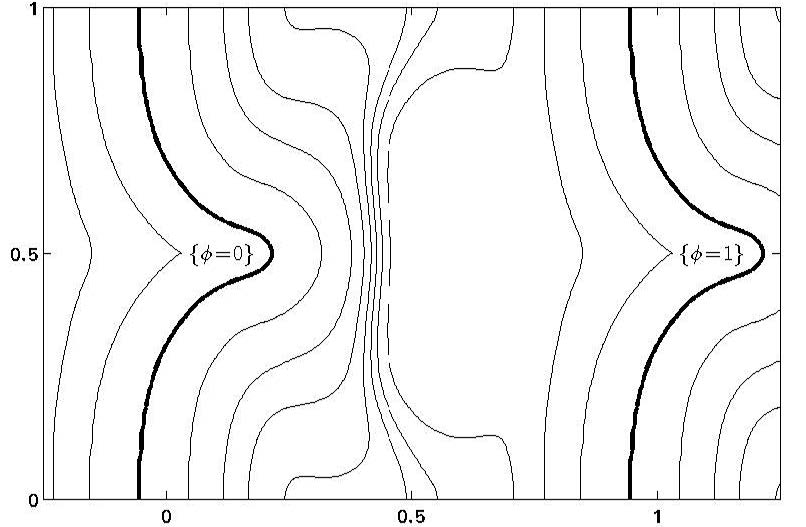}
\caption{(a) Graphs of $\bar{S}(\phi)$ and $\bar{D}(\phi)$. (b) Contour plot of the testing $\phi(x,t)$. (c) Reinitialized  $\phi(x,t)$ without smoothing at $t=0.2$. (d) Reinitialized  $\phi(x,t)$ at $t=0.2$ with one smoothing iteration every time step. (e) Reinitialized  $\phi(x,t)$ at $t=0.2$ with 10 smoothing iterations every time step.}
\label{Reini}
\end{figure}

When the flame front travels very fast, the level set function becomes very flat. When the motion of the flame front nearly stops, the level set function becomes very sharp. 
In either case the computational error will increase, and the level set may not be well captured. 
Hence reinitialization needs to be applied regularly to keep the level set function approximately equal to the signed distance function near the level set.

The standard reinitialization equation is
\be\label{Re1}
\phi_t+S(\phi)(|D\phi|-1)=0\ ,\ S(\phi)=\mathrm{sgn}(\phi),
\ee
which spreads out the signed distance from the level set $\{\phi(x,t)=0\}$. 
The function $S:\Rm\ra\Rm$ can be mollified to improve the numerical accuracy, see chapter 7 of \cite{OF02} for details.

To perform reinitialization on (\ref{Tu1}) with $P=e_1$, $\phi(x,t)$ must satisfy $\phi(x+e_1,t)=\phi(x,t)+1$ and be 
periodic in $x_2,\dots,x_n$. 
See Figure \ref{Reini}.(b) for an example of $\phi(x,t)$ using contour plot. 
To maintain the spatial periodicity, we modify (\ref{Re1}) as follows:
\be\label{Re2}
\phi_t+\bar{S}(\phi)(|D\phi|-1)=0,
\ee
where $\bar{S}$ is a 1-periodic function and $\bar{S}(\phi)=\mathrm{sgn}(\phi)$ for $\phi\in[-1/2,1/2]$. 
See Figure \ref{Reini}.(a) for the graph of the mollified version of $\bar{S}(\phi)$. 
In numerical computation, (\ref{Re2}) is discretized by WENO5 scheme and RK3 schemes with time step $\Delta t=\Delta x$. 

However, Figure \ref{Reini}(c) shows that (\ref{Re2}) spreads out distances from both $\{\phi=0\}$ and $\{\phi=1\}$. 
As a result, $\phi(x,t)$ is squeezed near $\{\phi=1/2\}$. 
The computation grinds to a halt when $\phi(x,t)$ becomes too sharp. 
To avoid this problem, we consider the following nonlinear diffusion equation:
\be\label{Re3}
\phi_t=c\bar{D}(\phi)\Delta\phi,
\ee
where $c$ is some positive constant and $\bar{D}:\Rm\ra\Rm$ is a 1-periodic function satisfying $\bar{D}(\phi)=0$ for $\phi\in [-\eps,\eps]$ and $\bar{D}(\phi)=1$ for $\phi\in [2\eps,1-2\eps]$. 
See Figure \ref{Reini}.(a) for the graph of $\bar{D}(\phi)$. 
Equation (\ref{Re3}) smooths $\phi(x,t)$ in the region away from the level set. 
In summary, we combine (\ref{Re2}) and (\ref{Re3}) to obtain the following reinitialization equation for the transformed G-equation (\ref{Tu1}) with $P=e_1$:
\be\label{Re4}
\phi_t+\bar{S}(\phi)(|D\phi|-1)=c\bar{D}(\phi)\Delta\phi.
\ee

In actual computation, we do not solve (\ref{Re4}) accurately because the diffusion term reduces the time step to $\Delta t=O((\Delta x)^2)$. 
Instead, we alternate between (\ref{Re3}) and (\ref{Re2}). Approximate (\ref{Re3}) by the simple iteration:
\be\label{Re5}
\phi_{i,j}:=(1-\bar{D}(\phi_{i,j}))\phi_{i,j}
+\bar{D}(\phi_{i,j})
{(\phi_{i+1,j}+\phi_{i-1,j}+\phi_{i,j+1}+\phi_{i,j-1})\over 4}.
\ee
The iteration (\ref{Re5}) is repeated a few times in each time step of the numerical scheme of (\ref{Re2}). 
This way, the time step remains $\Delta t=O(\Delta x)$. 
See Figure \ref{Reini}(d) and (e) for an illustration of the smoothing effect.

\section{Numerical Results}
\setcounter{equation}{0}
We consider all G-equations (\ref{Gi}),(\ref{Gc}),(\ref{Gv}),(\ref{Gs}) in two spatial dimensions with $P=e_1$ and $s_L=1$. 
The velocity field $V(x,t)$ is chosen to be cellular flow (\ref{Cell}) with various values of the intensity $A$ to study the growth rate of turbulent flame speed. 
Also the Markstein number $d$ is varied to study the curvature and strain effect.

First we solve the periodic initial value problem (\ref{Tu1}) for $G(x,t)$ on $[0,1]^2$. 
Then by (\ref{Tu2}) we construct the solution $G(x,t)$ in some stripe domain $[a,b]\times [0,1]$ and obtain the level set $\{G(x,t)=0\}$. 
The computational domain is $[0,1]\times[0,1]$ with grid points up to $400\times 400$. 

\begin{figure}
\center (a) Inviscid G-equation  \\
\includegraphics[width=0.9\textwidth]{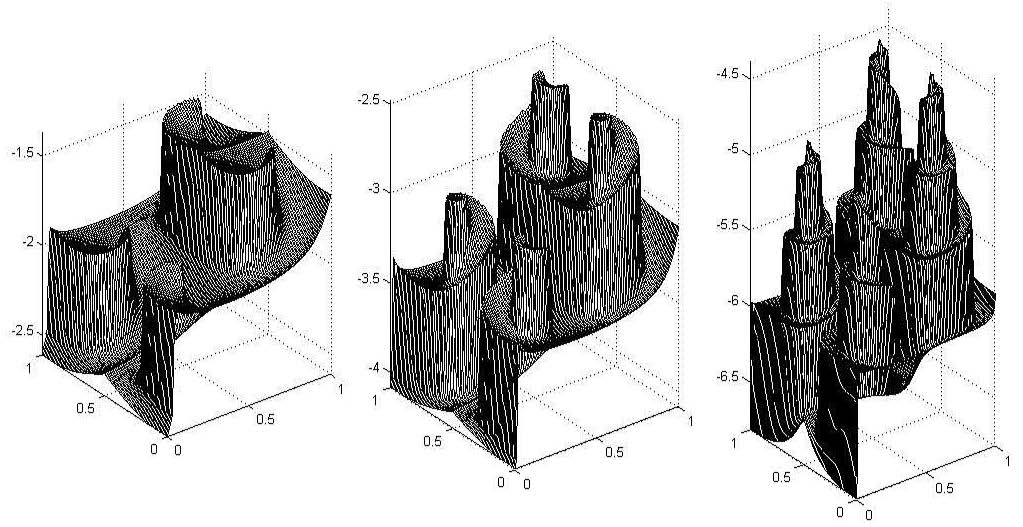}
\center (b) Curvature G-equation \\
\includegraphics[width=0.9\textwidth]{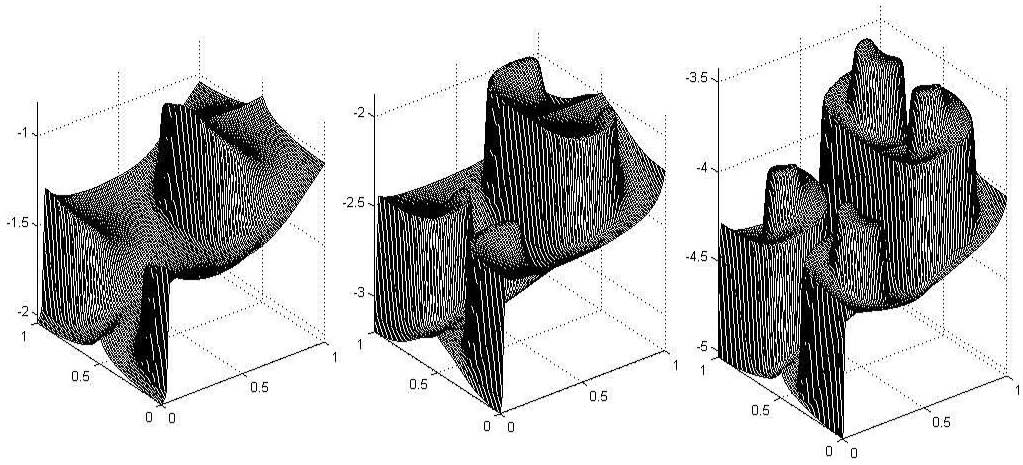}
\center (c) Viscous G-equation   \\
\includegraphics[width=0.9\textwidth]{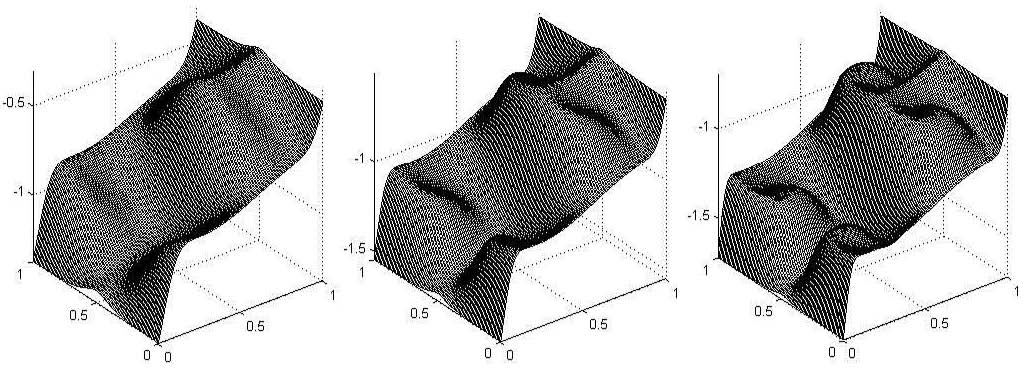}
\caption{Graphs of $G(x,t)$ at $t=1$ for inviscid, curvature, and viscous G-equations in 
cellular flow with $A=4,8,16$ (left to right) and $d=0.1$.}
\label{Graph}
\end{figure}

Figure \ref{Graph} shows the graphs of $G(x,t)$ for inviscid, curvature, and viscous G-equations at $t=1$ with $A=4,8,16$ and $d=0.1$. 
When $A$ is large, the graph of $G(x,t)$ has cone shape in each cell. 
Due to the curvature effect, $G(x,t)$ is less irregular and the cone formation is slower. 
The regular viscosity makes $G(x,t)$ even smoother.

\begin{figure}
\center (a) Inviscid G-equation  \\
\includegraphics[width=0.8\textwidth]{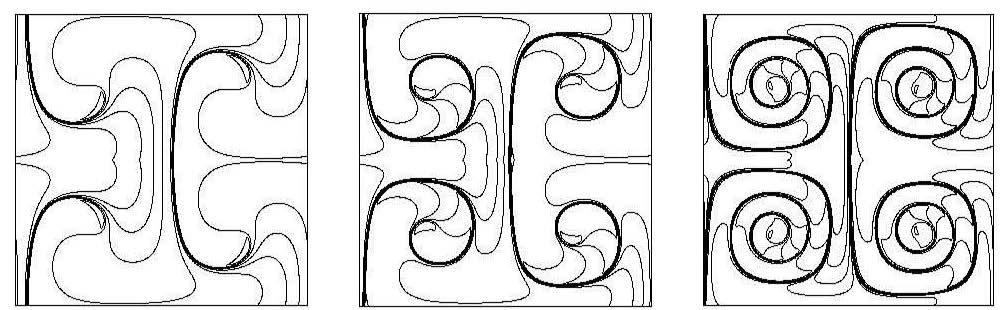}
\center (b) Curvature G-equation \\
\includegraphics[width=0.8\textwidth]{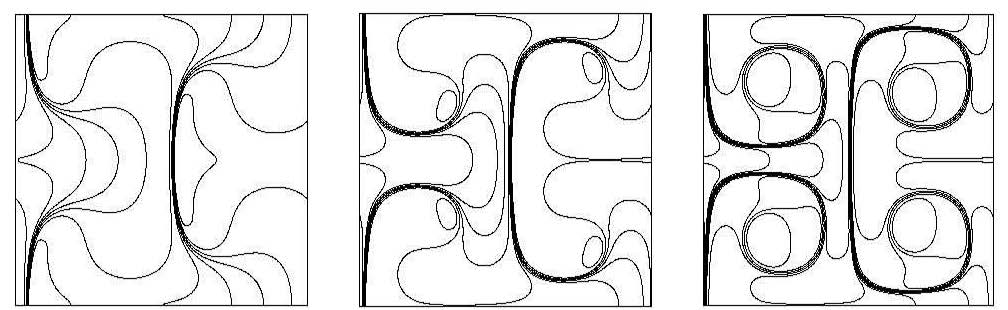}
\caption{Contour plots of $G(x,t)$ at $t=1$ for inviscid and curvature G-equations in 
cellular flow with $A=4,8,16$ (from left to right) and $d=0.1$.}
\label{Contour}
\end{figure}

Figure \ref{Contour} shows the contour plots of $G(x,t)$ for inviscid and curvature G-equations. 
When the level set merges, shock waves occur and the derivative of $G(x,t)$ is discontinuous across the shock wave. We observe that the shock wave is of spiral shape in each cell, especially at $d=0.1$, $A=16$.

\begin{figure}
\center (a) Inviscid G-equation \hspace{0.75in}
(b) Curvature G-equation \\
\includegraphics[width=0.48\textwidth]{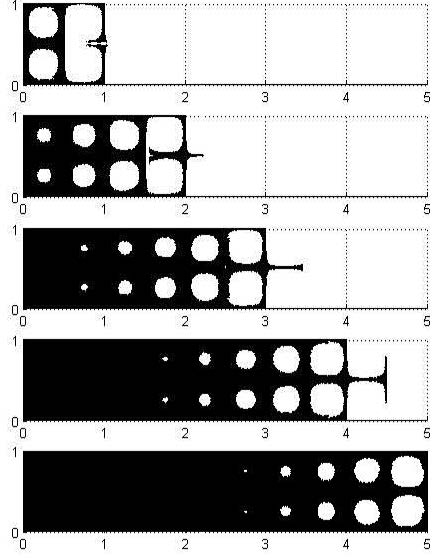}
\includegraphics[width=0.48\textwidth]{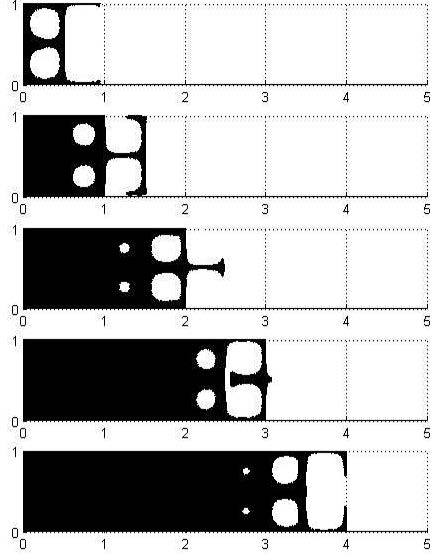}
\caption{Propagation of flame front in time for inviscid and curvature G-equations with $A=32$, $d=0.1$ and $t= 0.1, 0.2, 0.3, 0.4, 0.5$.} 
\label{Flow}
\end{figure}

Figure \ref{Flow} shows the propagation of the flame front for inviscid and curvature G-equations at $A=32$ and $d=0.1$. 
When $A$ is large, the flame front of the inviscid G-equation travels faster along the boundaries of the cells with bubbles formed behind. 
The flame front spirals inside the cells, and the bubbles shrink in the wake. 
If the curvature effect is added, the flame front is concave when traveling along the boundaries. 
The curvature term slows down front propagation yet the wake bubbles shrink faster.

\begin{figure}
\center (a) Inviscid G-equation  
(b) Curvature G-equation (c) Viscous G-equation\\
\includegraphics[width=\textwidth]{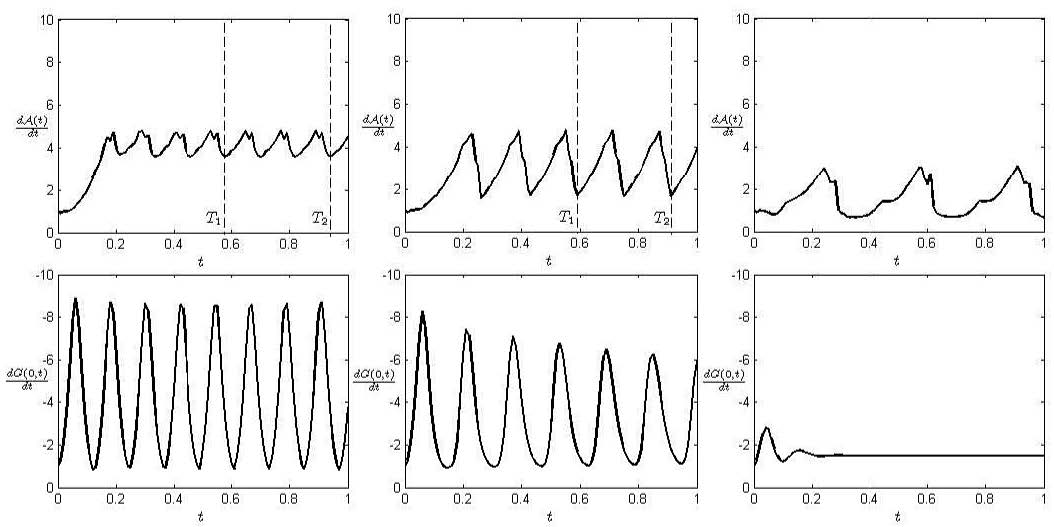}
\caption{Plots of $\mathcal{A}'(t)$ and $G'(0,t)$ for inviscid, curvature and viscous G-equation with $A=8$ and $d=0.1$.} 
\label{Hbar}
\end{figure}

Figure \ref{Hbar} shows the time derivative function of $\mathcal{A}(t)$ and $G(x=0,t)$ for inviscid, curvature and viscous G-equation with $A=8$ and $d=0.1$. 
After a short time interval, $\mathcal{A}'(t)$ behaves like a periodic function.
Hence we can approximate $s_T$ by taking the average of $\mathcal{A}'(t)$ over a periodic time interval:
$$
s_T \approx {1\over T_2-T_1}\int_{T_1}^{T_2}\mathcal{A}'(t)dt
={\mathcal{A}(T_2)-\mathcal{A}(T_1)\over T_2-T_1}.
$$ 
See Figure \ref{Hbar} for examples of selections of $T_1$, $T_2$. 
So we don't need to use (\ref{Tu5}) and perform large time simulation in order to approximate $s_T$ correctly.

Next we consider the behavior of $G'(x,t)$ in time for fixed $x$. ($'$: $\partial/\partial t$.)
For inviscid G-equation, $G'(x,t)$ behaves like a periodic function after a short time, hence we can evaluate $s_T$ by the same method as above rather than using (\ref{Tu8}). 
For viscous G-equation, the dissipation term causes damping in $G'(x,t)$. 
Hence  $G'(x,t)$ converges to $-s_T$ in time, and $G(x,t)$ converges to the stationary solution (\ref{Tu7}). 
For curvature G-equation, however, we see only slight damping in $G'(x,t)$. 

\begin{figure}
\center (a) Inviscid G-equation \hspace{0.4in}
(b) Curvature G-equation \\
\includegraphics[width=0.8\textwidth]{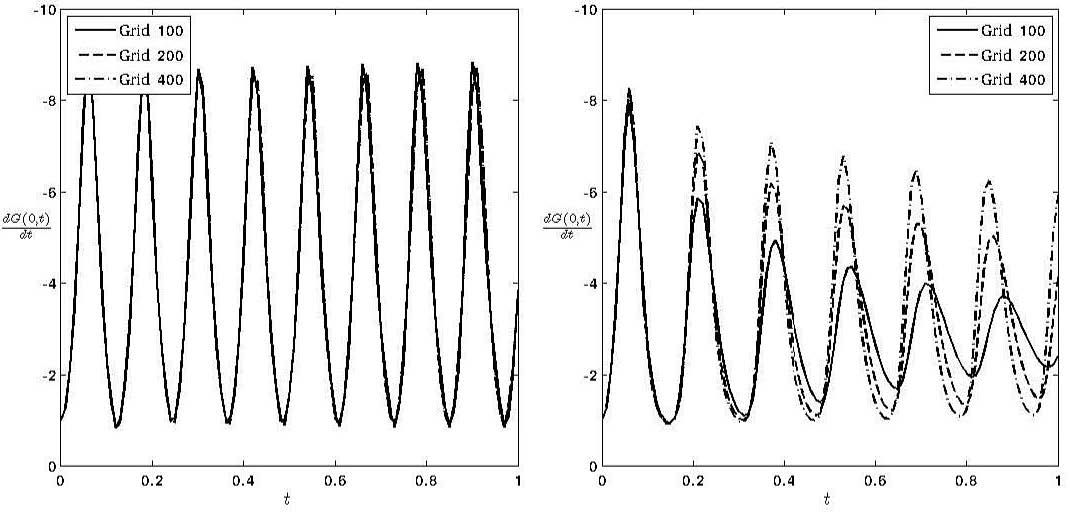}
\caption{Plots of $G'(0,t)$ for inviscid and curvature G-equation with $A=8$, $d=0.1$ and grid sizes 100, 200, 400.} 
\label{HDamp}
\end{figure}

Figure \ref{HDamp} shows function $G'(0,t)$ with different grid sizes. 
For inviscid G-equation, the numerical scheme is higher order accurate, and the artificial dissipation is well minimized. 
Hence damping is hardly observed even on coarse grid. 
For curvature G-equation, the numerical scheme is second order accurate, and the curvature term may be incorrectly evaluated at shock wave. 
Hence damping effect is very strong on coarse grid, and we must use fine grid to reduce the artificial diffusion.

We denote $s_{T}^{inv}$, $s_{T}^{cur}$, $s_{T}^{vis}$, $s_{T}^{str}$ the turbulent flame speeds for inviscid, curvature, viscous, strain G-equations respectively. 
We also denote them as functions of either the flow intensity $(A)$ or the Markstein number $(d)$. 
Note that when $A=0$ we have
$s_{T}^{inv}=s_{T}^{cur}(d)=s_{T}^{vis}(d)=s_{T}^{str}(d)=s_L$,
and when $d=0$ we have
$s_{T}^{inv}(A)=s_{T}^{cur}(A)=s_{T}^{vis}(A)=s_{T}^{str}(A)$.

\begin{figure}
\center (a) \hspace{2in}
(b)\\
\includegraphics[width=0.48\textwidth]{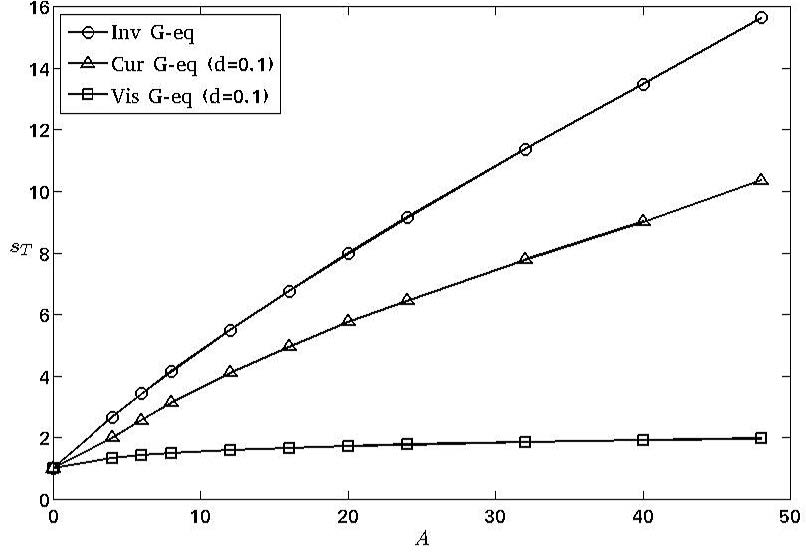}
\includegraphics[width=0.48\textwidth]{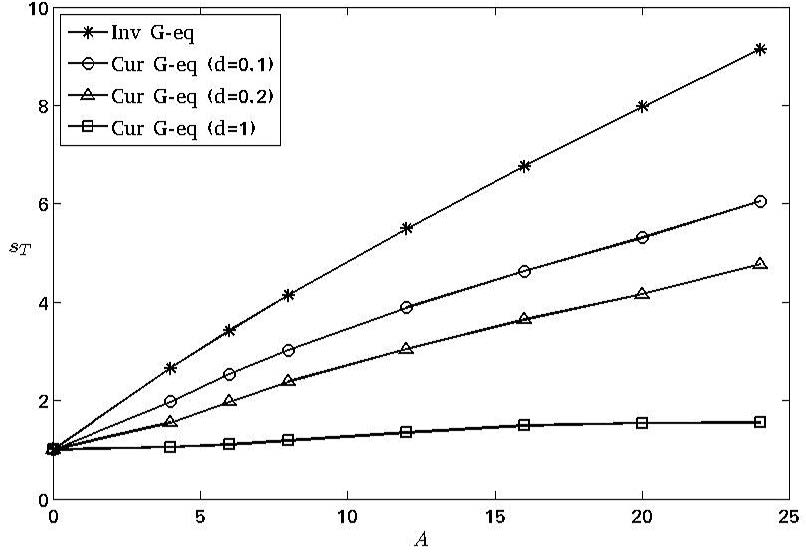}
\caption{(a) Plots of $s_T=s_T(A)$ for inviscid, curvature and viscous G-equations with $d=0.1$. (b) Plots of $s_T=s_T(A)$ for curvature G-equation with $d=$ 0.1, 0.2 and 1.}
\label{InCuVi}
\end{figure}

Figure \ref{InCuVi}(a) shows the graphs of $s_{T}^{inv}(A)$, $s_{T}^{cur}(A)$ and $s_{T}^{vis}(A)$ with $d=0.1$. 
The numerical results indicate that they all increase as $A$ increases and
$$
s_{T}^{vis}(A)\leq s_{T}^{cur}(A)\leq s_{T}^{inv}(A).
$$
Figure \ref{InCuVi}(b) shows the graphs of $s_{T}^{inv}(A)$ and $s_{T}^{cur}(A)$ with $d=0.1,0.2,1$. 
We used the forward Euler scheme for $d=$ 0.1 and semi-implicit scheme for $d=$ 0.2 and 1. 
It is known that $s_{T}^{inv}(A)=O({A/\log A})$ and $s_{T}^{vis}(A)=O(1)$.
However, the precise asymptotic behavior of $s_{T}^{cur}(A)$ as $A\ra \infty$ remains open. 
The growth scaling of $s_{T}^{cur}(A)$ is not conclusive from the range of $A$ we simulated. 

\begin{figure}
\center (a) Strain G-equation ($d=0.01$) \hspace{0.5in}
(b) Strain G-equation ($d=0.02$) \\
\includegraphics[width=0.48\textwidth]{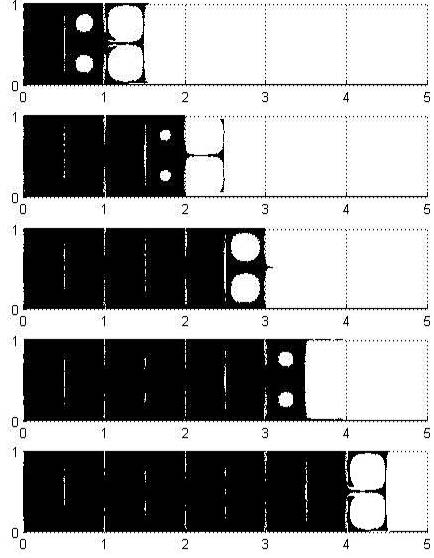}
\includegraphics[width=0.48\textwidth]{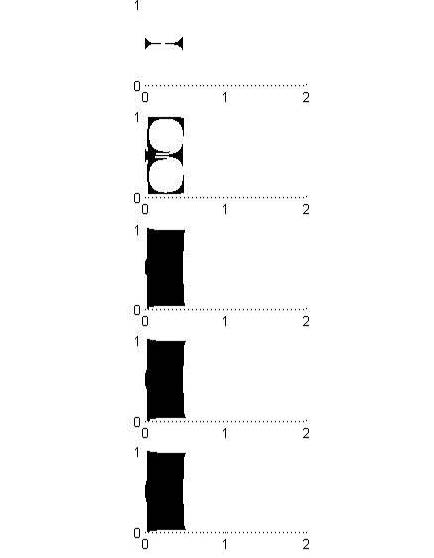}
\caption{Propagation of flame front in time for strain G-equation with $A=32$, $d=0.01,0.02$ and $t=0.3,0.6,0.9,1.2,1.5$.} 
\label{FlowStr}
\end{figure}

Figure \ref{FlowStr} shows the propagation of the flame front for strain G-equation with $A=32$ and $d=$ 0.01, 0.02. 
Near the corner of the cell, the velocity field is weak ($|V(x)|\approx 0$) yet the strain rate is strong ($|\mathcal{S}|\approx 2\pi A$). 
In the strong advection scheme, the strain term dominates near the corner of the cell, and the flame front cannot reach the corner. 
At $d=0.01$, Figure \ref{FlowStr}(a) shows that incomplete combustion occurs near the corners of the cells, yet the flame front still manages to propagate. 
At $d=0.02$, however, the flame front stops moving after $t=0.6$. 
Note that if the level set stops moving, then the level set function forms a sharp layer. 
Here reinitialization is needed to alleviate the stiff level set function and keeps computation going.

\begin{figure}
\center (a) \hspace{2in} (b)\\
\includegraphics[width=0.48\textwidth]{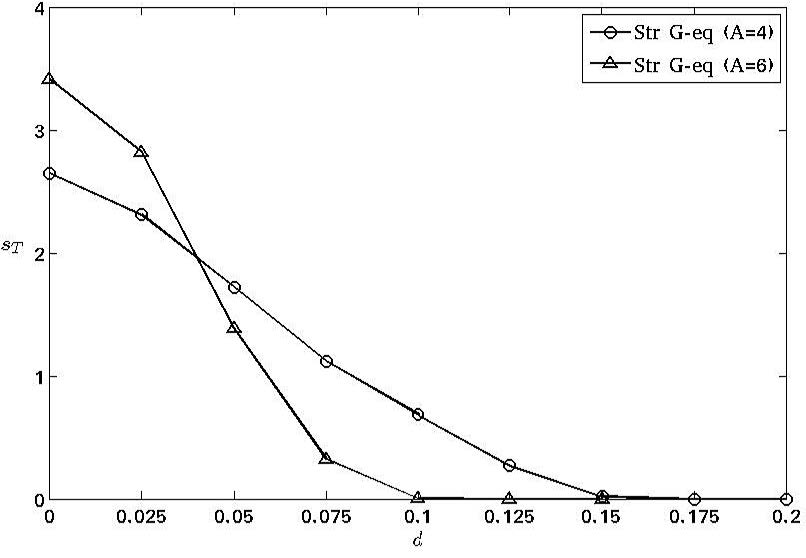}
\includegraphics[width=0.48\textwidth]{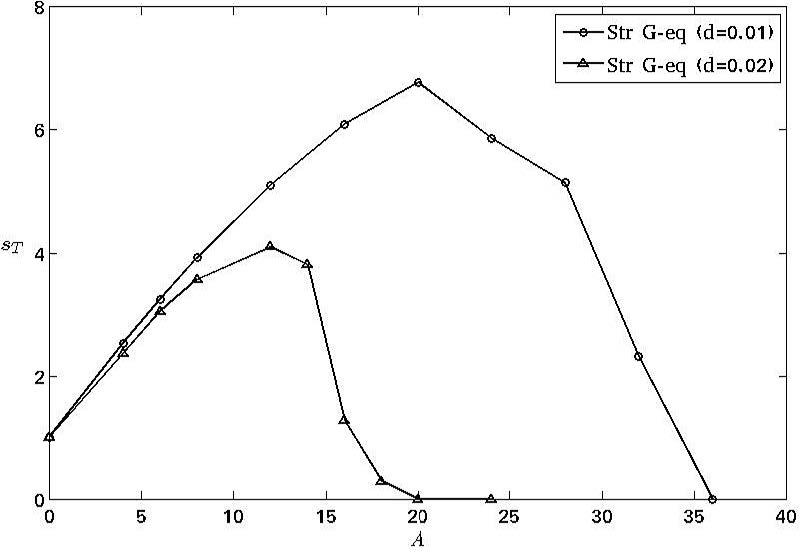}
\caption{(a) Plots of $s_T=s_T(d)$ for strain G-equation with $A=4,6$. (b) Plots of $s_T=s_T(A)$ for strain G-equation with $d=0.01,0.02$.}
\label{Strain}
\end{figure}

Figure \ref{Strain}(a) shows the graphs of $s_{T}^{str}(d)$ with $A$= 4, 6. 
In contrast to $s_{T}^{vis}(d)\geq s_L$ for any $d>0$ \cite{LXY10}, $s_{T}^{str}(d)$ decreases to zero when $d$ is large enough. Figure \ref{Strain}(b) shows the graphs of $s_{T}^{str}(A)$ with $d=$ 0.01, 0.02. 
When $A$ is small, $(s_L-d\mathcal{S})$ remains positive and $s_{T}^{str}$ is increasing. 
When $A$ gets larger, $s_{T}^{str}$ decreases and eventually drops down to zero. 
This agrees with the nonlinear phenomenon in turbulent combustion that high strain is the cause of flame quenching \cite{B92, R95}.

\section{Conclusion}
\setcounter{equation}{0}
We have studied various G-equation models numerically, and 
evaluated the corresponding turbulent flame speeds in cellular flows. 
Based on the numerical results, we showed how the turbulent flame speeds are affected by viscosity, curvature or strain effect. Weak and strong bending effects of the speeds caused are observed in curvature and viscous G-equations.
Quenching effect only appears in the strain G-equation. 
In future work, we plan to study turbulent flame speeds of G-equations in time dependent or three dimensional spatially periodic vortical flows.

\appendix
\renewcommand{\theequation}{\Alph{section}.\arabic{equation}}

\section{Appendix: Surface Stretch Rate Formula}
\setcounter{equation}{0}

\begin{figure}
\center \includegraphics{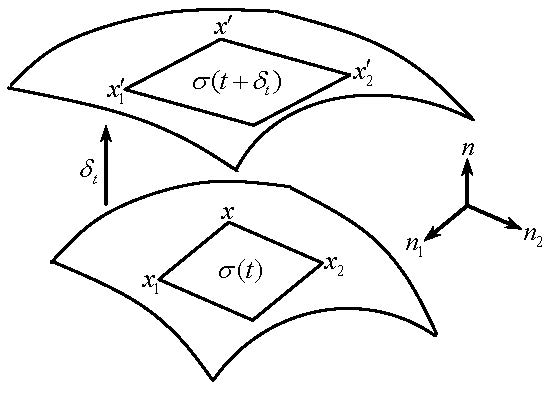}
\caption{An illustration of surface stretch in the proof of Appendix A.} 
\label{Surface}
\end{figure}

In this appendix, we derive the surface stretch rate. 
A surface stretch rate formula in three dimensions is derived in \cite{M83}. 
Here we give an alternative formula in any dimensions and apply it in G-equation.

\begin{theo}
Suppose a smooth hypersurface in $\Rd$ is moving in the velocity field $V(x,t)$. 
Denote $\sigma$ the surface element area and $n$ the unit normal vector of a point on the surface. 
Then the surface stretch rate is given by
\be\label{Su1}
{1\over\sigma}{d\sigma\over dt}
=\mathrm{div}(V)-n\cdot DV\cdot n.
\ee
\end{theo}

\nit {\bf Proof}: 
See figure \ref{Surface} for the picture of the proof. 
Fix a time $t$ and a point $x$ on the surface, the surface can be locally approximated by its tangent plane. Let $\{n_1,\dots,n_{d-1}\}$ be an orthonormal basis of the tangent plane and $\eps_1,\dots,\eps_{d-1}$ be infinitesimal scalars. 
Then the surface element can be presented by a rectangle whose sides are the vectors $\eps_1n_1,\dots,\eps_{d-1}n_{d-1}$. 
The surface element area is
$$
\sigma(t)=\eps_1\cdots\eps_{d-1}.
$$
For $1\leq k\leq d-1$, denote $x_k=x+\eps_kn_k$ the neighboring point of $x$ of the rectangle. 
Then we say the rectangle is determined by the staring point $x$ and neighboring points $x_1,\dots,x_{d-1}$. 

After a time interval $\dta_t$, suppose the new locations of $x$, $x_k$ are $x'$, $x'_k$ respectively. 
Then the surface element becomes a parallelogram determined by the staring point $x'$ and neighboring points $x'_1,\dots,x'_{d-1}$. 
Denote $\dta_k=x_k'-x'$, then the sides of the parallelogram are the vectors $\dta_1,\dots,\dta_{d-1}$. 
The surface element area is
$$
\sigma(t+\dta_t)=\sqrt{\det(A^TA)},
$$
where $A=[\dta_1,\cdots,\dta_{d-1}]$ is the matrix whose columns are the sides of the parallelogram.

From now on we keep all calculations up to first order of $\dta_t$ and omit higher order terms. 
The surface moves in velocity field $V(x,t)$, then
$$
x'=x+V(x,t)\dta_t\ ,\ x_k'=x_k+V(x_k,t)\dta_t.
$$
$$
\Rightarrow \dta_k=(x_k-x)+DV\cdot(x_k-x)\dta_t
=(\mathbb{I}_d+\dta_tDV)\cdot (\eps_k n_k).
$$
Denote $\eta_k=(\mathbb{I}_d+\dta_tDV)\cdot n_k$ and $N=[n_1,\cdots,n_{d-1}]$, then $\dta_k=\eps_k\eta_k$ and
$$
\sigma(t+\dta_t)=\eps_1\cdots\eps_{d-1}\sqrt{\det(B^TB)},
$$
where $B=[\eta_1,\cdots,\eta_{d-1}]=(\mathbb{I}_d+\dta_tDV)N$.
Then we have
$$
B^TB=N^T(\mathbb{I}_d+\dta_tDV^T)(\mathbb{I}_d+\dta_tDV)N=\mathbb{I}_{d-1}+\dta_tN^T(DV+DV^T)N
$$
$$
\Rightarrow \det(B^TB)=1+\dta_t\mathrm{tr}(N^T(DV+DV^T)N)
=1+2\dta_t\mathrm{tr}(N^TDVN)
$$
$$
\Rightarrow \sigma(t+\dta_t)=\sigma(t)(1+\dta_t\mathrm{tr}(N^TDVN)).
$$

Hence the surface stretch rate is 
$$
\lim_{\dta_t\ra0}{1\over\sigma(t)}{\sigma(t+\dta_t)-\sigma(t)\over \dta_t}
=\mathrm{tr}(N^TDVN))
$$
$$
=n_1^TDVn_1+\cdots+n^T_{d-1}DVn_{d-1}.
$$
Note that $\{n_1,\dots,n_{d-1},n\}$ is an orthonormal basis of $\Rd$, then
$$
n^T_1DVn_1+\cdots+n^T_{d-1}DVn_{d-1}+n^TDVn
=\mathrm{tr}(DV)=\mathrm{div}(V).
$$
We combine the last two equations and finish the proof. \qed

\begin{rmk}
Result of \cite{M83} in three dimensions reads:
\be\label{Su2}
{1\over\sigma}{d\sigma\over dt}
=(n\cdot V)\mathrm{div}(n)-\mathrm{curl}(V\times n)\cdot n.
\ee
Indeed we can verify that (\ref{Su1}) and (\ref{Su2}) are equivalent in $\Rm^3$.
\end{rmk}

\begin{cor}
Let $V(x,t)$ be an incompressible flow and denote $\kappa=\mathrm{div}(n)$ the curvature of the surface.
If the surface moves in the velocity field $V(x,t)$ and the normal direction with constant speed $s_L$:
\be\label{Su3}
{dx\over dt}=V(x,t)+s_Ln,
\ee
then the stretch rate is
\be\label{Su4}
{1\over\sigma}{d\sigma\over dt}
=-n\cdot DV\cdot n+s_L\kappa.
\ee
\end{cor}
\nit {\bf Proof}: Substitute (\ref{Su3}) into (\ref{Su1}), then we have
$$
{1\over\sigma}{d\sigma\over dt}
=\mathrm{div}(V)+s_L\mathrm{div}(n)
-n\cdot DV\cdot n
-s_Ln\cdot Dn\cdot n.
$$
The first term is $0$ due to incompressibility of $V$. 
By some calculations, the last term is $0$. \qed

\section{Appendix: Iteration Scheme for Cell Problem of Viscous G-equation}
\setcounter{equation}{0}

In this appendix, we prove the convergence of the iteration scheme for the cell (corrector) problem of viscous G-equation at large enough $d$: 
$$
-ds_L\Delta u^{(k+1)}+V(x)\!\cdot\! Du^{(k+1)}
= H^{(k)}-s_L|P+Du^{(k)}|-V(x)\!\cdot\! P,\ x\in\Tset^n,
$$
\be\label{It1}
H^{(k)}=s_L\int_{\Tset^n}|P+Du^{(k)}|dx.
\ee

First we verify the solvability of (\ref{It1}). 
Denote $L^2_{per}$ and $H^1_{per}$ the spaces of all mean zero and periodic functions in $L^2(\Tset^n)$ and $H^1(\Tset^n)$ respectively. 
Since $V(x)$ is assumed to be periodic, mean zero and divergence free, by Fredholm alternative theorem, the equation
$$
-\Delta u+V(x)\cdot Du=f\ ,\ x\in\Tset^n
$$
has unique weak solution $u\in H^1_{per}$ provided $f\in L^2_{per}$. 
If $u^{(k)}\in H^1_{per}$ then the right hand side of (\ref{It1}) is in $L^2_{per}$ and there exists unique solution $u^{(k+1)}\in H^1_{per}$ for (\ref{It1}). 
Therefore given any $u^{(1)} \in H^1_{per}$ then we can construct a sequence $\{u^{(k)}\}_{k\in\Nset}$ in $H^1_{per}$.

\begin{theo}
The sequence $\{u^{(k)}\}_{k\in\Nset}$ in $H^1_{per}$ defined by the iteration scheme (\ref{It1}) converges provided $d>\sqrt{n}/\pi$. 
\end{theo}

\nit {\bf Proof}: Replace the index $k$ in (\ref{It1}) by $k+1$ and take their difference:
$$
-ds_L\Delta (u^{(k+2)}-u^{(k+1)})+V(x)\cdot D(u^{(k+2)}-u^{(k+1)})
$$
$$
=(H^{(k+1)}-H^{(k)})-s_L\left[|P+Du^{(k+1)}-|P+Du^{(k)}|\right].
$$
Multiply the equation by $u^{(k+2)}-u^{(k+1)}$ and take integration over $\Tset^n$:
$$
d\int_{\Tset^n} \left[D(u^{(k+2)}-u^{(k+1)})\right]^2dx
$$
\be\label{It2}
= -\int_{\Tset^n} \left[|P+Du^{(k+1)}|-|P+Du^{(k)}|\right]
(u^{(k+2)}-u^{(k+1)})dx.
\ee
Here we use the fact that $V(x)$ is divergence free and $u^{(k+2)}-u^{(k+1)}$ is mean zero.
Recall the Poincar\'e inequality:
$$
\norm{u}_{L^2(\Tset^n)}\leq {\sqrt{n}\over\pi} \norm{Du}_{L^2(\Tset^n)}
\ , \ u\in H^1_{per}.
$$
By Cauchy inequality, (\ref{It2}) implies that
$$\begin{array}{rl}\displaystyle
&d\norm{Du^{(k+2)}-Du^{(k+1)}}^2_{L^2(\Tset^n)}
\vspace{0.05in}\\
\leq& \displaystyle
\norm{|P+Du^{(k+1)}|-|P+Du^{(k)}|}_{L^2(\Tset^n)}
\norm{u^{(k+2)}-u^{(k+1)}}_{L^2(\Tset^n)}
\vspace{0.05in}\\
\leq& \displaystyle
\norm{Du^{(k+1)}-Du^{(k)}}_{L^2(\Tset^n)}
{\sqrt{n}\over\pi}
\norm{Du^{(k+2)}-Du^{(k+1)}}_{L^2(\Tset^n)},
\end{array}
$$
$$
\Rightarrow \norm{Du^{(k+2)}-Du^{(k+1)}}^2_{L^2(\Tset^n)}
\leq {\sqrt{n}\over\pi d}
\norm{Du^{(k+1)}-Du^{(k)}}_{L^2(\Tset^n)}.
$$
If $d>\sqrt{n}/\pi$, then $\{Du^{(k)}\}_{k\in\Nset}$ is contracting in $L^2(\Tset^n)$. By Poincar\'e inequality, $\{u^{(k)}\}_{k\in\Nset}$ converges in $H^1_{per}$. \qed

\bibliographystyle{plain}

\end{document}